\documentclass[12pt]{amsart}
\usepackage{amssymb,amsmath}
\numberwithin{equation}{section}
\def\simleq{\underset\sim<}

\def\T{\text}
\def\1#1{\overline{#1}}
\def\2#1{\widetilde{#1}}
\def\3#1{\widehat{#1}}
\def\4#1{\mathbb{#1}}
\def\5#1{\frak{#1}}
\def\6#1{{\mathcal{#1}}}
\def\C{{\4C}}
\def\R{{\4R}}

\def\Z{{\4Z}}

\begin{document}
%
\title[separate real analyticity and CR extendibility]
{separate real analyticity and CR extendibility}
\author[ L.~Baracco G.~Zampieri ]
{Luca Baracco Giuseppe Zampieri }
\maketitle
\def\B{\Bbb B}
\def\Giialpha{\mathcal G^{i,i\alpha}}
\def\cn{{\C^n}}
\def\cnn{{\C^{n'}}}
\def\ocn{\2{\C^n}}
\def\ocnn{\2{\C^{n'}}}
\def\const{{\rm const}}
\def\rk{{\rm rank\,}}
\def\id{{\sf id}}
\def\aut{{\sf aut}}
\def\Aut{{\sf Aut}}
\def\CR{{\rm CR}}
\def\GL{{\sf GL}}
\def\Re{{\sf Re}\,}
\def\Im{{\sf Im}\,}
\def\codim{{\rm codim}}
\def\crd{\dim_{{\rm CR}}}
\def\crc{{\rm codim_{CR}}}
\def\phi{\varphi}
\def\eps{\varepsilon}
\def\d{\partial}
\def\a{\alpha}
\def\b{\beta}
\def\g{\gamma}
\def\G{\Gamma}
\def\D{\Delta}
\def\Om{\Omega}
\def\k{\kappa}
\def\l{\lambda}
\def\L{\Lambda}
\def\z{{\bar z}}
\def\w{{\bar w}}
\def\Z{{\1Z}}
\def\t{{\tau}}
\def\th{\theta}
\emergencystretch15pt
\frenchspacing
\newtheorem{Thm}{Theorem}[section]
\newtheorem{Cor}[Thm]{Corollary}
\newtheorem{Pro}[Thm]{Proposition}
\newtheorem{Lem}[Thm]{Lemma}
\theoremstyle{definition}\newtheorem{Def}[Thm]{Definition}
\theoremstyle{remark}
\newtheorem{Rem}[Thm]{Remark}
\newtheorem{Exa}[Thm]{Example}
\newtheorem{Exs}[Thm]{Examples}
\def\Label#1{\label{#1}}
\def\bl{\begin{Lem}}
\def\el{\end{Lem}}
\def\bp{\begin{Pro}}
\def\ep{\end{Pro}}
\def\bt{\begin{Thm}}
\def\et{\end{Thm}}
\def\bc{\begin{Cor}}
\def\ec{\end{Cor}}
\def\bd{\begin{Def}}
\def\ed{\end{Def}}
\def\br{\begin{Rem}}
\def\er{\end{Rem}}
\def\be{\begin{Exa}}
\def\ee{\end{Exa}}
\def\bpf{\begin{proof}}
\def\epf{\end{proof}}
\def\ben{\begin{enumerate}}
\def\een{\end{enumerate}}

\def\simto{\overset\sim\to\to}
\def\1alpha{[\frac1\alpha]}
\def\T{\text}
\def\R{{\Bbb R}}
\def\I{{\Bbb I}}
\def\C{{\Bbb C}}
\def\Z{{\Bbb Z}}
\def\Fialpha{{\mathcal F^{i,\alpha}}}
\def\Fiialpha{{\mathcal F^{i,i\alpha}}}
\def\Figamma{{\mathcal F^{i,\gamma}}}
\def\Real{\Re}

\section{Introduction}
In $\C^2=\R^2+i\R^2$  with coordinates $z=(z_1,z_2),\,\,z=x+iy$, we consider a function $f$ continuous on a domain $\Omega$  of $\R^2$ 
separately real analytic in $x_1$ and CR  extendible to $y_2$ (resp. CR extendible to $y_2>0$). This means that  $f(\cdot,x_2)$ extends holomorphically for $|y_1|<\epsilon_{x_2}$ and $f(x_1,\cdot)$ for $| y_2|<\epsilon$ (resp. $0\leq y_2<\epsilon$ continuous up to $y_2=0$) with $\epsilon$ independent of $x_1$. We prove  in Theorem~\ref{t2.2} that $f$ is then real analytic (resp.  in Theorem~\ref{t2.3} that it extends holomorphically to a ``wedge" $W= \Omega+i\Gamma_\epsilon$ where $\Gamma_\epsilon$ is an open cone  trumcated by $|y|<\epsilon$ and containing the ray $0<y_2<\epsilon)$. The extension is uniformly continuous up to $y=0$ and gives $f$ as  limit.
Note that the first result can be obtained from the second: it follows, via the edge of the wedge theorem of \cite{AH81} from the holomorphic extension to the pair of wedges $W^\pm$ which correspond to the two sides $\pm y_2>0$. But we prefer to give its own simplified proof.
We point out that it is not made any assumption on uniform continuity or uniform boudedness for the different extensions: these come as consequences of the statement.  In any case  boundedness is the 
main issue: once this is proved, then continuity follows (Theorem~\ref{t1.1}), and the holomorphic extension 
of a function which is CR extendible to  both  $y_1$ and $y_2$
is a consequence of the edge of the wedge theorem of Ajrapetyan-Henkin \cite{AH81}. Historically, it was earlier obtained by Malgrange-Zerner with specification of the side: separate CR continuous extension to $y_1>0$ and $y_2>0$ implies holomorphic extension to the quadrant $y_1>0,\,\,y_2>0$. (See also Komatsu \cite{K72} and Druzkowski \cite{D80} where continuity is replaced by boundedness.) 
Successively, Siciak \cite{S69}, \cite{S69bis} proved that CR extension in both arguments implies real analyticity without  assumption of continuity or boundedness. What we show in the present paper is that CR extension to $y_2$ and simple real analyticity in $x_1$ suffices. And that CR extension to $y_2>0$ implies holomorphic extension to a wedge of $\C^2$ with edge  $\R^2$ and which points to that side. We stress attention to the fact that CR extension in at least one argument is needed as shows the classical example of  $f(x_1,x_2)=x_1x_2\T{exp}(-\frac1{x_1^2+x_2^2})$.
As for the proof of our  result, the main technical tool is a refined version of the Hartogs Lemma which is contained in Theorem~\ref{hartogslemma}: it is a combination of the Fatou's Theorem and the Phragm\'en-Lindel\"of principle. 
In the other respects, to tackle the problem, we introduce a  point of view which is new in the context of the separate analyticity, that is the CR theory. 
First, we show that the CR extensions in either variable are continuous when the other argument is restricted to an open dense range of values. We then apply the edge of the wedge theorem to extend $f$ to an open set of $\C^2$ and next the theorem by Hanges-Treves \cite{HT81} to propagate this extension along the planes $y_2=\T{const}$. At this point we apply our generalized Hartogs Lemma to the sequence of the roots of the Taylor coefficients of the expansions of $f$ in $z_2$ to fill all missing values of $x_2$ in the domain of $f$.

The authors  are grateful to professor Claudio Rea for fruitful discussions. 
\newline
MSC: 32D10, 32U05, 32V25

\section{generalized/uniformly-continuous boundary values}
Let $\R^n$ and $\C^n$ be the euclidean real and complex spaces with coordinates $x$ and $z=x+iy$ respectively and 
also write $x=(x',x_n),\,\,z=(z',z_n)$. 
Let $\Omega$  be a domain of $\R^n$. We discuss of boundary values on $\Omega$ of functions which solve the Cauchy-Riemann equations on manifolds or wedges of $\C^n$ with boundary or edge $\Omega$. 
 We begin by a result, essentially due to \cite{R86},  on uniform continuity 
of functions which 
are holomorphic and tempered in one variable, say $z_n$, for $y_n>0$. 
\bt
\Label{t1.1}
Let $f$ be a continuous function on $\Omega$ such that 
for some $\epsilon$ and for all $x'$, $f(x',\cdot)$ extends holomorphically in $z_n$ to $0\leq y_n<\epsilon$, uniformly continuous up to $y_n=0$. We also suppose that for any $x^o_n$ there are two neighborhoods $A\subset\subset B$ of $x^o_n$ such that $|f(x',x_n+iy_n)|\simleq |y_n|^{-k}$ $\forall x',\,\,\forall x_n\in B\setminus A$ and for a suitable $k$. Then $f$ is uniformly continuous for $0\leq y_n\leq \epsilon'$ and $x\in\tilde\Omega\subset\subset\Omega$.
\et
\bpf
Since $f|_\Omega$ is already known to be continuous, it suffices to prove that
$$
f(x',x_n+iy_n)-f(x',x_n)\to0\T{ as }y_n\to0\T{ uniformly on any $\tilde\Omega\subset\subset\Omega$.}
$$
We take a function $\chi=\chi(z_n)$, of class $C^\infty_c$, $\chi\equiv1$ at $x^o_n$, such that $\T{supp}(\bar\partial\chi)\cap\R\subset B\setminus A$ and $\bar\partial\chi=O(|y_n|^k)$. Under this choice we have that $\bar\partial(\chi f)$ is uniformly bounded. 
We write $F(\tau):=\chi f(x',x_n+i\tau y_n)$ and apply Cauchy formula to the function $\frac{F(\tau)}{\tau+1}$ at $\tau=1$ for the half-plane $\Pi^+=\{\Real \tau>0\}$. We get, after substituting $\zeta=\xi+i\eta$ for $i\tau$,
\begin{equation}
\Label{1.1}
\chi f(x',x_n+iy_n)=\frac1\pi\int_{-\infty}^{+\infty}\frac{\chi f(x',x_n+\xi y_n)}{\xi^2+1}d\xi -\frac2\pi\int\int_{\Pi^+}-i\bar\partial_{z_n}(\chi f)(x',x_n+\zeta y_n) y_n\frac1{\zeta^2+1} d\xi d\eta.
\end{equation}
Recall that $\frac1\pi\int_{-\infty}^{+\infty}\frac1{\xi^2+1}d\xi=1$. Hence \eqref{1.1} implies
\begin{equation}
\Label{1.2}
\chi f(x',x_n+iy_n)-\chi f(x',x_n)=\frac1\pi\int_{-\infty}^{+\infty}\frac{\chi f(x',x_n+ \xi y_n)-\chi f(x',x_n)}{\xi^2+1}d\xi+O(|y_n|),
\end{equation}
which yields the conclusion.
\epf

Theorem~\ref{t1.1} can be extended to a general wedge swept by analytic discs. This is locally described as the image, through a CR embedding into $\C^n$, of $W=\Omega+i\Gamma_\epsilon$ where $\Gamma_\epsilon$ is an open convex cone trumcated by $|y|<\epsilon$ in a linear subspace of $\R^n$. 
Here $\epsilon=\epsilon_x$ has a positive lower bound on each compact subset of $\Omega$.
The embedding is supposed to be smooth up to $\Omega$. It also holds for more general sets swept by analytic half-discs.
The set $W=\Omega+i\Gamma_{\epsilon}$ is said a ``wedge" with ``edge" $\Omega$ and ``directional cone" $\Gamma$. 
Any subset $\tilde W=\tilde\Omega+i\tilde\Gamma_{\tilde\epsilon}$ with $\tilde\Omega\subset\subset\Omega$, $\tilde\Gamma\subset\subset\Gamma$ 
(that is $\overline{\tilde\Gamma}\setminus\{0\}\subset\Gamma$)
and $\tilde\epsilon<\epsilon$ is said a ``proper" subwedge of $W$; we will use the notation $\tilde W\subset\subset W$ in this case. 
When $\Gamma$ is open in $\R^n$,  Theorem~\ref{t1.1} can be remarkably improved.  
We want to explain in what sense a function, or more generally a distribution or a hyperfunction, $f$ on $\Omega$ is the boundary value of a holomorphic function that we still denote $f$ on a wedge $W$ of dimension $2n$. A generalized (or hyperfunction)
 boundary value acts on $g$ in the space $\mathcal A(\bar\Omega)$ of real analytic functions on $\bar\Omega$ through the integral
\begin{equation}
\Label{1.10}
\int_{\Omega}f(x+iy)g(x+iy)\,dx,
\end{equation}
 for $y\in\tilde\Gamma\subset\subset\Gamma$ suitably small. By Cauchy formula, the integrations for different $y$ only ``differ" near $\partial\Omega$.
 Thus \eqref{1.10} defines an element in the quotient $\frac{\mathcal A'(\bar\Omega)}{\mathcal A'(\partial\Omega)}$ that we denote by $\T{bv}(f)$.
  If the holomorphic function $f$ on $W$ is tempered in the sense that  $|f|\simleq |y|^{-k}$ for some $k$ as $y\to0, \,\,y\in\tilde\Gamma$, then it has a distributional boundary value defined over $\phi\in C^{k+1}_c(\Omega)$ by
\begin{equation}
\Label{1.11}
\langle \T{\rm bv}(f),\phi\rangle=\underset{y\to 0\,\,y\in\tilde\Gamma}{\T{\rm lim }}\int_{\Omega}f(x+iy)\phi(x )\,dx.
\end{equation}
In this case \eqref{1.11} is compatible with \eqref{1.10}
in the sense that the hyperfunction defined on $\Omega$ by \eqref{1.10} coincides with the distribution defined by \eqref{1.11}.
For continuous functions  on $\Omega$ we have coincidence of all possible notions of boundary values:
\bt
\Label{t1.2} (Rosay \cite{R86} and Cordaro - private communication)
Let $W=\Omega+i\Gamma_\epsilon$ be a wedge of  dimension $2n$ and let $f$ be a continuous function on $\Omega$ which is the boundary value of a holomorphic function on $W$ that we still denote by $f$. 
Then,
 $f$ is uniformly continuous on any proper subwedge $\tilde W\subset\subset W$ up to $\tilde\Omega$ and the boundary value is in fact a limit.
\et
\bpf
The direction in the wave front set of $f$ are contained in the polar cone $\Gamma^*$ and this is the same as a 
hyperfunction  and a distribution. In this situation $f$ is the boundary value of a tempered holomorphic function on 
$W$
(cf. \cite{H84} Theorem 8.4.15) which must coincide, 
by uniqueness, with the former extension. At this point one can apply 
the analogous of Theorem~\ref{t1.1} for functions which are holomorphic in a wedge 
(\cite{R86} Proposition 1, point 3).

\epf
Differently from ``boundary value" the word ``extension" is ambiguous. An ``extension" has full meaning when it is a holomorphic function on a  wedge $W$ of dimension $2n$ and reproduces the initial function $f$ through \eqref{1.10} or \eqref{1.11}. As we have just seen, when $f$ is continuous on $\Omega$ and  is the boundary value of a holomorphic function on a wedge $W$ of dimension $2n$, then it is in fact its uniform limit.
But a function may happen to extend along discs which fill a wedge $W$ 
of general dimension
without being the uniform limit of its ``extension". First, the separate extensions may not glue into a continuous function on $W$ as in the following example which 
was suggested to us by professor P. Cordaro.
\be
We take $\Omega=\R\times\R$ and define $f$ on $\Omega$ by $f(x_1,x_2)=x_1\,\T{\rm sin }\frac{x_2}{x_1}$. Then $f$ is continuous on $\Omega$ and extends as an entire function  along each plane $\{x_1\}\times \C$:
thus we have here, instead of a wedge, a manifold without boundary $\R\times\C$. 
 But it is not tempered in the neighborhood of any point $(0,x_2)$ as we may check by using for instance the sequence 
$$
(x_1^\nu,z_2^\nu)=(\frac1{\nu^2},x_2+i\frac1\nu).
$$
\ee
Also, even if  $f$ extends separately along discs which cover a wedge $W$ of dimension $2n$ and the different extensions glue into a holomorphic function on $W$, there is no evidence that $f$ is the boundary value, hence the uniform limit, of its extension. However, this conclusion holds, all over $\Omega$, when we know from the beginning that it is true over a sufficiently large subset of $\Omega$.  
\bt
\Label{t1.3}
In the situation of Theorem~\ref{t1.1}, assume in addition that $f$ extends holomorphically to a wedge $W=\Omega+i\Gamma_\epsilon$ for an open cone $\Gamma\subset\R^n$ which contains the positive $y_n$-axis. Then, $f$ is uniformly continuous on any $\tilde W\subset\subset W$ up to $\tilde\Omega$.
\et
\bpf
Since $f$ is holomorphic in $W$, 
we know in particular that it is continuous for $y_n>0$.
 In combination with Theorem~\ref{t1.1}, we get the uniform continuity for  
 $y'=0,\,\,y_n\geq0,\,\,x\in\tilde\Omega$.
In this situation, if we inspect \eqref{1.10}, we see that $f$ is the boundary value, in the sense of hyperfunctions, of its extension. To conclude we have just to apply Theorem~\ref{t1.2}.

\epf

\section{two-sided and one-sided holomorphic extension from $\R$ to $\C$}
Before entering the main theme, we need some preliminar result on functions of one complex variable.
Let $\C$ be the complex plane with coordinate $\tau$, $\Delta$  the unit disc $\Delta=\{\tau\in\C:\,|\tau|<1\}$, $\Delta^+$ the upper half-disc $\Delta^+=\{\tau\in\Delta:\,\Im \tau\geq0\}$, $I$ the unit real interval $I=\{t\in\R:\,|t|<1\}$. We recall that a subharmonic function is a function which is upper semicontinuous and whose value at any point $\tau$ is dominated by the mean value on the boundary of any disc centered at $\tau$. The following generalization of Hartogs Lemma will play a crucial role.
\bt
\Label{hartogslemma}
Let us be given a sequence $\{\phi_\nu\}_\nu$ of functions 
upper semicontinuous on $\bar\Delta^+$, subharmonic on $\overset\circ\Delta^+$ and which satisfy for some constants $l$ and $L$
\begin{equation}
\Label{2.1}
\begin{cases}
\underset{\nu}{\T{\rm lim~sup } }\underset{\tau\in\partial\Delta^+}{\T{\rm sup }}\phi_\nu(\tau)\leq l\T{  $\forall\tau\in\partial \Delta^+$},
\\
\underset{\nu}{\T{\rm lim~sup } }\phi_\nu(t)\leq0\,\,\forall t\in I,
\\
\underset{\tau\in\partial\Delta^+}{\T{\rm sup }}\phi_\nu(\tau)\leq L\,\,\forall\nu.
\end{cases}
\end{equation}
Then for any $\alpha$ and $\eta$ there is $\nu_{\alpha\,\eta}$ such that $\forall  \nu\geq\nu_{\alpha\,\eta}$
\begin{equation}
\Label{2.2}
\phi_\nu(\tau)\leq\alpha+l\kappa\,\Im \tau\quad\T{ for $\Im \tau\geq\eta$ }
\end{equation}
where $\kappa $ is an universal constant.
\et
\bpf
We denote by $\chi$ the function on $\partial\Delta^+$ which is $0$ for $\Im \tau=0$ and $l$ for $|\tau|=1$.
Let $P_z(\zeta),\,\,z\in\Delta^+,\,\zeta\in\partial\Delta^+$, be the Poisson kernel of $\Delta^+$. 
For any $z\in\overset\circ\Delta^+$ and $\alpha>0$ and for suitable $\nu_{\alpha,z}$ we have when $\nu\geq \nu_{\alpha\,z}$:
\begin{equation}
\begin{split}
\phi_\nu(z)&
\leq \int_{\partial\Delta^+}P_z(\zeta)\phi_\nu(\zeta)d\lambda(\zeta)
\\
&\leq \int_{\partial\Delta^+}P_z(\zeta)\chi(\zeta)d\lambda(\zeta)+\alpha
\\
&\leq \kappa l\,\Im z+\alpha,
\end{split}
\end{equation}
where the first inequality comes from the subharmonicity of each $\phi_\nu$, the second from Fatou's Lemma and the third from the Phragm\'en-Lindel\"of principle.
For any $\alpha$, $\eta$ and for suitable $\gamma_{\alpha\,\eta}$ we have
\begin{equation}
\Label{poisson}
\left|\int_{\partial\Delta^+}(P_z(\zeta)-P_w(\zeta))\phi_\nu(\zeta)d\zeta\right|<\alpha\T{ for $|w-z|<\gamma_{\alpha\,\eta}, |\Im z|\geq\frac\eta2,\,\,|\Im w|\geq \frac\eta2$},
\end{equation}
because $\{P_z(\zeta)\}_{\zeta\in \partial\Delta^+}$, as a family of functions of $z$, is equicontinuous on $|\Im z|>\frac\eta2$ 
and since the $\phi_\nu$'s are uniformly bounded.
By a finite covering of $\Delta^+\cap\{\Im \tau\geq\eta\}$ by discs $\Delta_{\gamma_\alpha}$
of radius $\gamma_\alpha$ and center at points $z$ with  $|\Im z|\geq\frac\eta2$, we get the conclusion of the proof of the Theorem.

\epf
\br
 We will get the same conclusion \eqref{2.2}  if we replace the half-disc  $\Delta^+$ by the strip $U^+_\delta=I+i(0,\delta)$.
\er

 Let $z=x+iy$ be the coordinates in $\C^n$, $\Omega$  an open domain of $\R^n$ and $f$ a continuous function on $\Omega$. 
\bd
\begin{itemize} We adopt the following terminology.
\item
We say that $f$ is separately real analytic in $x_j$ if its restriction to the section of $\Omega$ with each line parallel to the $x_j$-axis is real analytic. Thus when all the other coordinates are fixed, $f$ extends to $|y_j|<\epsilon_x$.
\item
We say that $f$ is separately CR extendible to $y_j$ if it is separately real analytic in $x_j$ and moreover it has holomorphic extension to $|y_j|<\epsilon$ with $\epsilon$ having positive lower bound locally uniform in $x$.
\item
We say that $f$ is separately CR extendible to $y_j>0$ or $y_j<0$ if it has a holomorphic extension to $0<y_j<\epsilon$ or $-\epsilon<y_j<0$ with $\epsilon$ locally uniform in $x$ and the extension on each plane $z_j$ has locally uniform limit $f$ for $y_j=0$. 
\end{itemize}
\ed
Let $\C^n=\C^{n_1}\times\C^{n_2}$ with coordinates $z=(z',z'')=(x'+iy',x''+iy'')$.
Separate analyticity in the group   of variables $x''$ means analyticity in $x''$ when the other variables are fixed.  
CR extendibility to all the directions of the $y''$-plane $\R^{n_2}$ means holomorphic extendibility to $|y''|<\epsilon$ with $\epsilon$ locally uniform in $x$.
CR extendibility to a cone $\Gamma_2\subset\R^{n_2}$  means that  $f$ extends holomorphically  for $y''\in\Gamma_2,\,\,|y''|<\epsilon$ with $\epsilon$ locally uniform in $x$; we also assume that, for fixed $x'$, $f$ is locally uniformly continuous up to $y''=0$. 
Clearly real analyticity or CR extendibility in a group of variables $x''$ is more restrictive than 
the combination of real analyticities
in any single $x''_j$.
We state now our main results whose proofs will follow in Section~4.
\bt
\Label{t2.2}
Let $f$ be a continuous function on $\Omega$ which is separately real analytic in $x'$ and CR extendible to $y''$. 
Then $f$ is real analytic.
\et
Theorem~\ref{t2.2} improves Siciak's theorem of \cite{S69} where it is assumed that 
$f$ is separately CR extendible both to $y'$ and $y''$.
By iteration, Theorem~\ref{t2.2} implies that if $f$ is separately real analytic in $x_1$ and CR extendible to each single variable $y_{2}$, ...,$y_n$, then it is in fact real analytic.
Let $\Omega=\Omega_1\times\Omega_2 \subset\R^n=\R^{n_1}\times\R^{n_2}$.
We will denote by $\Gamma$ an open convex cone in the $y$-space $\R^n$ and by $\Gamma_\epsilon$ its trumcature by $|
y|<\epsilon$; we will denote by $\Gamma_2$ or $\Gamma_{2\,\epsilon}$ the analogous open cones in the $y''$-space $\R^{n_2}$. 
We fix a cone $\Gamma_2$.
\bt
\Label{t2.3}
Let $f$ be a continuous function on $\Omega$ 
which is real analytic in  $x'$ and CR extendible to $\Gamma_2$. Then
 $f$ extends holomorphically to 
 a wedge $W=\Omega+i\Gamma_\epsilon$ where $\Gamma$ is a conic neighborhood of an arbitrarily large proper subcone of $\Gamma_2$ and is uniformly continuous  up to $y=0$.
\et
\bc
\Label{c2.1}
Let $\Omega_1\times V_2$ be a domain in $\R^{n_1}\times \C^{n_2}$ and let $f$ be a  function in $\Omega_1\times V_2$, 
such that  $\forall x'\in\Omega_1$, $f(x',\cdot)$ is
holomorphic and $\forall z''\in V_2$, $f(\cdot,z'')$ is real analytic. Then $f$ extends holomorphically to a neighborhood $W$ of $\Omega_1\times V_2$.
\ec
(Cf. Theorem~1 by Shiffman \cite{Sh89}.) 
This is a corollary to Theorem~\ref{t2.2} or Theorem~\ref{t2.3} applied to the family of domains $\Omega$ obtained by slicing $\Omega_1\times V_2$ by the planes $y''=\T{const}$:  some care is needed because $f|_{\Omega_1\times V_2}$ is not assumed to be continuous.  When $V_2$ is a wedge $\Omega_2+i\Gamma_{2\,\epsilon}$, one can see that $W$ is also a wedge $\Omega+i\Gamma_\epsilon$ 
where $\Gamma\supset\tilde\Gamma_2$ for an arbitrarily large $\tilde\Gamma_2\subset\subset\Gamma_2$.
 However, Theorem~\ref{t2.3} is far better: real analyticity in $x'$ for fixed $x''\in\Omega_2$, and not  $z''\in \Omega_2+i\Gamma_{2\,\epsilon}$, suffices.

Let us notice that Theorem~\ref{t2.2} could be obtained  from  Theorem~\ref{t2.3}. In fact, CR extension to $y''$ implies CR extension to any pair of antipodal open cones $\Gamma^\pm_2\subset\R^{n_2}$. Application of 
Theorem~\ref{t2.3}
implies extension to a pair of wedges $W^\pm$, whose profiles are conical neighborhoods $\Gamma^\pm$ of $\Gamma^\pm_2$
in $\R^{n}$ . But then the Ajrapetyan-Henkin edge of the wedge theorem  gives extension to the directions of the convex hull of the profiles that is the whole of the directions of the $y$-plane $\R^n$. However, for the sake of clearness, we will give a separate proof.

\noindent
{\bf Outline of the proof of Theorems~\ref{t2.2}, \ref{t2.3}.}
 We start from Theorem~\ref{t2.2}. 
Assuming that $\Omega$ is the square $I^2=(-1,1)\times(-1,1)$ of $\R^2$ and denoting by $\Delta$ the standard disc,  we consider the 
sets $\underset{x_1\in I}\cup\{x_1\}\times \Delta$ and $\underset{x_2\in I}\cup\Delta_{\epsilon_{x_2}}\times\{x_2\}$ 
 to which $f$ is supposed to extend. 
 Now, 
 by Baire's theorem  the extension of $f$ is  bounded in $I_\delta\times \Delta$ and $U_\delta\times I_\delta$ for some interval $I_\delta=(-\delta,\delta)$  and some strip $U_\delta=I+iI_\delta$. 
 It is easy to prove, by Cauchy's inequalities,
  that boundedness combined with continuity on $I^2$ implies in fact continuity for the extensions.
  By the celebrated Ajrapetyan-Henkin's ``edge of the wedge" Theorem of \cite{AH81}, $f$ extends to $\Delta_\delta\times\Delta_\delta$ for a new $\delta$.
 The argument of this theorem consists in ``attaching" analytic discs to the two hypersurfaces and in extending $f$ along these discs by the maximum principle.
 Also, wherever the extension is continuous, we can apply the propagation of the holomorphic extendibility of CR functions along  discs: in particular $f$ will extend holomorphically to $U_\delta\times\Delta_\delta$ for a new $\delta$.  
 At this point we apply our generalized Hartogs Lemma 
 - Theorem~\ref{hartogslemma} above -
 to the sequence 
 $\{\phi_\nu(z_1)\}_\nu=\{\left|\frac{\partial_{z_2}^\nu(z_1,0)}{\nu!}\right|^{\frac1\nu}\}_\nu$ 
 and get normal convergence of the Taylor series of $f$ in $z_2$ over discs of radius arbitrarily close to $1$ and uniform for $|y_1|$ suitably small. Thus $f$ extends holomorphically to $U_\epsilon\times\Delta$ (where $\Delta$ is a little shrunk).
 
 As for Theorem~\ref{t2.3}, where we have extension only to the side $y_2>0$, we follow the same lines. The proof that boundedness implies continuity requires an extra argument since we are now working in the upper half-disc $\Delta^+$ at the boundary points $\Im z_2=0$. This is provided by 
 the more general Theorem~\ref{t1.1}; a  direct proof, based on the
 Phragm\'en-Lindel\"of priciple, could be given. In this way we get holomorphic extension to $U_\delta\times\Delta^+_\delta$. We center now the Taylor series of $f$ with respect to $z_2$ at  
 $z_2=\frac\delta2-\sigma+i\frac\delta2$  for $\sigma<<\delta$ 
 and get that this series converges normally for 
 $z_2$ in the disc $\Delta_{\delta'}(\frac\delta2-\sigma+i\frac\delta2)$ of center $\frac\delta2-\sigma+i\frac\delta2$ and radius $\delta'<\frac\delta2$ and
  uniformly for $z_1$ belonging to compact subsets of $I+i(I_{\epsilon_\delta}\setminus\{0\})$ for a suitable $\epsilon_\delta$. 
  We next iterate and center the Taylor series of (the extension of) $f$ to the points $j\delta'-j\sigma$ or $-j\delta'+j\sigma\,\,j=0,...,N$ so that $|N\delta'-N\sigma|>1$. In this way we are able to extend $f$ to the domains $W_\delta:=(I+iI_{\epsilon^N_\delta})\times V_\delta$ where $V_\delta$ is a neighborhood of $I+i\frac\delta2$ in the $z_2$-plane. By taking the union $\cup_\delta W_\delta$ we obtain a domain of extension 
 which is a neighborhood of $I^2+i(\{0\}\times I^+_{\delta})$. Next, by a theorem by Kashiwara (cf. \cite{K72}), we get extension to a domain $W=I^2+i\Gamma_\epsilon$ where $\Gamma$ is a conical neighborhood of the half-line $y_2>0$ in $\R^2$.
\section{Proofs}
\noindent
{\bf Proof of Theorem~\ref{t2.2}.}
By using iteration, we may assume $n=2$. 
 We will use the notations $I=(-1,1)$, $I^2=I\times I$, $I_\epsilon=(-\epsilon,\epsilon)$, $\Delta=\{\tau: |\tau|<1\}$, $\Delta_\epsilon=\{\tau:|\tau|<\epsilon\}$,
 $U_\epsilon=\{\tau\in\C:\,\,|\Real \tau|<1,\,\,|\Im \tau|<\epsilon\},$
$\Delta^+=\{\tau\in\Delta:\,\Im \tau\geq0\}$, $I^+=\{t\in I:\,t\geq 0\}$, $U^+_\epsilon=\{\tau\in U_\epsilon:\,\,\Im\tau>0\}$. 
 The statement being local we can suppose that $f$ is defined and continuous in $I^2$, extends to
 $\Delta_{\epsilon_{x_2}}\times\{x_2\}$ and $\{x_1\}\times\Delta$ holomorphic with respect to $\tau\in\Delta$ and $\Delta_{\epsilon_{x_1}}$ respectively. We even assume that $f$ extends indeed to discs of radius slightly bigger than $1$ or $\epsilon_{x_2}$.  
 
 \noindent
 {\bf (a)}
  We first prove  that there are an open interval $I_\delta\subset I$  
  and an open strip $U_\delta=I+iI_\delta$ 
  such that $f$ is continuous in $I_\delta\times\Delta$ and in $U_\delta\times I_\delta$ (cf. \cite{H73} and \cite{S69}). Hence it is a continuous CR function therein. 
We start from the proof of the continuity on $I_\delta\times\Delta$. 
Let $K_l=\{x_1\in I: \underset{z_2\in\Delta}{\T{sup}}|f(x_1,z_2)|\leq l\}$. We note that $K_l\subset K_{l+1}$ and that $\cup_lK_l=I$ since, for each $x_1$, $\underset{z_2\in\Delta}{\T{sup}}|f(x_1,z_2)|<+\infty$. We claim that
\begin{equation}
\Label{2.1}
\begin{cases}
K_l\T{ is closed},
\\
f\T{ is continuous on $ K_l\times\Delta$}.
\end{cases}
\end{equation}
In fact, let $x^\nu_1\to x^o_1$ with $x^\nu_1\in K_l$; we want to show that then $x^o_1\in K_l$. We use the notation
$F_\nu(z_2):=f(x_1^\nu,z_2)-f(x^o_1,z_2)$. The sequence $\{F_\nu\}_\nu$ is equicontinuous in a neighborhood of $\Delta$: in fact, remember that that $f$ was supposed to be holomorphic for $|z_2|$ slightly bigger than $1$: so the conclusion follows from the hypothesis of boundedness in addition to the Cauchy inequalities. 
We claim that $F_\nu\to0$. Otherwise, by the equicontinuity, there is a subsequence $\{F_{\nu_k}\}_k$ which converges to a limit $F\neq0$. But this limit is holomorphic in $\Delta$ and $0$ in $I$, a contradiction.
 This proves the claim and thus  \eqref{2.1} follows.
We can see now that by Baire's Theorem the union $\cup_lK_l$ being the whole $I$, the sets $K_l$ must contain an open interval for large $l$. 
Also, such an interval can be found in a neighborhood of any point. It needs not to contain $0$ but 
we may assume it, by means of a small translation, for the purpose of our proof.
Thus we can assume that $f$ extends as a continuous function on  $ I_\delta\times \Delta$ holomorphic in $z_2$: hence it is a continuous CR function therein.

We pass now to prove that $f$ is a continuous CR function on $U_\delta\times I_\delta$. For this purpose, we define $J_l=\{x_2:\,\,f(\cdot,x_2)\T{\rm extends to $|y_1|<\frac1l$ and $|f(\cdot,x_2)|<l\}$}$. 
In fact, if $x_2^\nu\to x^o_2$ with $x_2^\nu\in J_l$, then by boundedness, there is a subsequence which converges to a holomorphic function on $U_{\frac1l}$; this must be $f(\cdot,x_2^o)$. As before we have $f(\cdot,x_2)|\leq l$ and $f|_{U_{\frac1l}\times I_l}$ is continuous.
 By Baire's theorem we still conclude that for large $l$, the set $J_l$ contains an open interval that we can suppose to be centered at $0$. This concludes the proof of the claim.

 \noindent
  {\bf (b)}
  At this point we apply the Ajrapetyan-Henkin edge of the wedge Theorem and conclude that $f$ extends holomorphically to a domain of type $\Delta_\delta\times\Delta_\delta$.
 Since this is a crucial point here, we give the outline of the proof which follows \cite{T98}.
 We show first how to extend $f$ for $0\leq \Im z_1<\delta,\,\,0\leq\Im z_2<\delta$.
  In fact, choose smooth functions $y_j(e^{i\theta})\geq0$ with supp$(y_1)\subset[0,\pi]$, supp$(y_2)\subset[\pi,2\pi]$ and with unit mean value, take $(\lambda_j)$ with $0\leq\lambda_j<\delta\,\,\forall j$, write $y_\lambda=(\lambda_1y_1,\lambda_2y_2)$ and consider the discs $A_{x_o,\lambda}(\tau)$ which are the holomorphic extensions 
  of $(x_o-T_0y_\lambda)+iy_\lambda$ from $\tau=e^{i\theta}\in\partial\Delta$ to $\tau\in\Delta$. (Here $T_0$ is the Hilbert transform normalized by the condition $T_0(\cdot)(0)=0$.) Note that the boundaries of these discs, 
   corresponding to the values $\tau=e^{i\theta}$ of the parameter,
   are contained in the union of $\Delta^+\times I_\delta$ and $I_\delta\times\Delta^+$. 
   Also, the set of their centers $\{A_{x_o,\lambda}(0)\}=\{x_o+i\lambda\}$ is the  set described by $0\leq\Im z_1<\delta,\,\,0\leq\Im z_2<\delta$. On the other hand $f$ is uniformly approximated over the set 
   of the boundaries 
   by a sequence 
   of polynomials according to the Baouendi-Treves approximation theorem ( cf. Theorem 1 ch. 13 of \cite{B91}). This sequence  is also convergent in the inside of these discs, in particular in the set of their centers, by the
   maximum principle. The limit of the sequence provides the desired extension of $f$ to the first quadrant $0\leq \Im z_1<\delta,\,\,0\leq\Im z_2<\delta$; in the same way we prove extension to 
   the other quadrants 
    and conclude the proof of our claim.

\noindent
 {\bf (c)}
 We notice now that
\begin{itemize}
\item
$f$ is continuous and CR on $U_\delta\times I_\delta$,
\item
$U_\delta\times I_\delta$ is foliated by the complex leaves $\Sigma_{x_2}:=U_\delta\times\{x_2\}$ for $x_2\in I_\delta$,
\item $f$ extends to $\Delta_\delta\times\Delta_\delta$,
\item each leaf $\Sigma_{x_2}$ intersects $\Delta_\delta\times I_\delta$.
\end{itemize}
But then the propagation of the holomorphic extendibility of CR functions along complex leaves yields extension of $f$  to an open domain $U_\delta\times\Delta_\delta$ of $\C^2$ for small $\delta$. We notice here that $\Delta_\delta\times\Delta_\delta$ is swept by discs with boundary in the region where $f$ is bounded and $U_\delta\times\Delta_\delta$ by discs with boundary in the union of $U_\delta\times I_\delta$ and $\Delta_\delta\times\Delta_\delta$ where $f$ is also bounded. Hence, by maximum modulus principle, $f$ is bounded in $U_\delta\times \Delta_\delta$. Since by such an intervals $I_\delta$, where $f$ has (different) bounds, we cover an open dense set $D\subset I$, we conclude that $f$ is continuous up to $y=0$ over $I\times D$. This remark will be crucial in the proof of the subsequent Theorem~\ref{t2.3}. 

\noindent
{\bf (d)} We consider now the Taylor series of $f$ with respect to $z_2$ centered at $z_2=0$
\begin{equation}
\Label{2.8}
\sum_\nu\frac{\partial_{z_2}^\nu f(z_1,0)}{\nu!}z_2^\nu.
\end{equation}
This represents a holomorphic function on $U_\delta\times \Delta_\delta$; when $x_1$ is fixed in $I$ this extends holomorphically for $z_2\in \Delta$.
We write $\phi_\nu(z_1):=\frac1\nu\T{log}\,\frac{|\partial_{z_2}^\nu f(z_1,0)|}{\nu!}$ and note that these are subharmonic functions of $z_1$. We have 
\begin{equation}
\Label{2.10}
\underset\nu{\T{\rm lim sup}}\,\underset{z_1\in{U_\delta}}{\T{sup}}\phi_\nu(z_1)\leq -\T{log}\,\delta,
\end{equation} 
and hence in particular $\underset{z_1\in{U_\delta}}{\T{sup}}\phi_\nu(z_1)\leq L$ for some constant $L$ and for any $\nu$.  We also have
\begin{equation}
\Label{2.11}
\underset\nu{\T{lim sup }} \phi_\nu(x_1)\leq0\T{ for $x_1\in I$}.
\end{equation}
(In \eqref{2.10} and \eqref{2.11} the domains $\Delta$ and $I$ should be arbitrarily little shrunk.) But then Theorem~\ref{hartogslemma} applies  to the sequence of the $\phi_\nu$'s over the pair of half-strips $x_1+U^\pm_\delta$. It implies  that for any $\alpha>0$ and $\eta>0$ and for suitable $\nu_{\alpha\,\eta}$, we have
\begin{equation*}
\underset { x_1\in I,\, y_1\geq \eta }{\T{sup}}\phi_\nu(z_1)\leq\frac\alpha2-\delta^{-1}\kappa\T{log}\,\delta|y_1|\quad\forall \nu\geq \nu_{\alpha\,\eta}.
\end{equation*}
We have  in other words 
\begin{equation*}
\frac{\left|\partial_{z_2}^\nu f(z_1,0)\right|^{\frac1\nu}}{\nu!}\leq e^{\frac\alpha2}\delta^{-\kappa|y_1|\delta^{-1}}.
\end{equation*}

Let $\epsilon=\epsilon_{\alpha\,\delta}$ satisfy $e^{\frac\alpha2}\delta^{-\kappa\epsilon\delta^{-1}}<(1-\alpha)^{-1}$; then, the series \eqref{2.8} converges uniformly for $z_1$ on compact subsets of $\dot U_\epsilon=I+i(I_\epsilon\setminus\{0\})$, normally for $z_2$ in $\Delta_{1-\alpha}$ and its sum is therefore a holomorphic function on $\dot U_\epsilon\times\Delta_{1-\alpha}$. Since we already know that this function extends to $y_1=0$ when $z_2\in\Delta_\delta$, then it is in fact holomorphic on $U_\epsilon\times\Delta$.
The proof is complete.
\newline
{\bf Proof of Theorem~\ref{t2.3}.}
It is not restrictive to assume $n=2$, $\Omega=I\times I$ and $\Gamma_2=I^+$.
We may suppose that $f$ is a continuous function 
 on $I^2$ which extends holomorphically to $\{x_1\}\times\Delta^+$ and $\Delta_{\epsilon_{x_2}}\times\{x_2\}$ $\forall x_1$ and $x_2$ in $I$: we prove that it extends holomorphically to a domain  $I^2+i\Gamma_\epsilon$ where $\Gamma$ is an open cone of $\R^2$ around the positive $y_2$-axis.

\noindent{\bf (e)}
The first part of the proof follows the lines of Theorem~\ref{t2.2}. We begin by noticing that $f$ extends continuously to $U_\delta\times I_\delta$ and to $I_\delta\times \Delta^+$ for some $\delta$: the former is identical as in Theorem~\ref{t2.2}, the latter is a consequence of Theorem~\ref{t1.1}. We then apply the edge of the wedge  and the propagation theorems and conclude that $f$ extends holomorphically to $U_\delta\times \Delta^+_\delta$. 
We note here that $\Delta^+_\delta$ contains $\Delta_{\frac\delta2}(i\frac\delta2)$ the disc with center $i\frac\delta2$ and radius $\frac\delta2$. We will  write $\delta$ instead of $\frac\delta2$ in the following, and therefore suppose that $f$ extends to $\Delta\times\Delta_\delta(i\delta)$.

\noindent
{\bf (f)} We consider now the Taylor series of $f$ with respect to $z_2$ centered at the point $\delta-\sigma+i\delta$ for $\sigma<<\delta$:
\begin{equation}
\Label{2.13}
\sum_\nu\frac{\partial_{z_2}^\nu f(z_1,\delta-\sigma+i\delta)}{\nu!}z_2^\nu.
\end{equation}
This represents a holomorphic function on $\Delta\times \Delta_\sigma(\delta-\sigma+i\delta)$; when $x_1$ is fixed in $I$ this extends holomorphically for $z_2\in \Delta_\delta(\delta-\sigma+i\delta)$ that is the radius of convergence in $z_2$ increases from $\sigma$ to $\delta$. 
We write $\phi_\nu(z_1):=\frac1\nu\T{log}\,\frac{|\partial_{z_2}^\nu f(z_1,\delta-\sigma+i\delta)|}{\nu!}$ and note that these are subharmonic functions of $z_1$. We have 
\begin{equation}
\Label{2.14}
\underset\nu{\T{\rm lim sup }}\underset{z_1\in{U_\epsilon}}{\T{sup}}\phi_\nu(z_1)\leq -\T{log}\,\sigma
\end{equation} 
together with
\begin{equation}
\Label{2.15}
\underset\nu{\T{lim sup }} \phi_\nu(x_1)< -\T{log}\,\delta \T{ for $x_1\in I$}.
\end{equation}
(In \eqref{2.14} and \eqref{2.15} the domains $U_\delta$ and $I$ should be arbitrarily little shrunk.) But then Theorem~\ref{hartogslemma} applies  to the sequence of the $\phi_\nu$'s over the pair of half-strips $U_\delta^\pm$. It implies  that for any $\alpha>0$ and $\eta>0$ and for suitable $\nu_{\alpha\,\eta}$, we have
\begin{equation*}
\underset { x_1\in I,\, y_1\geq \eta }{\T{sup}}\phi_\nu(z_1)\leq -\T{log}\,\delta-\kappa\T{log}\,\sigma|y_1|+\frac\alpha2\quad\forall \nu\geq \nu_{\alpha\,\eta}.
\end{equation*}
We have  in other words 
\begin{equation}
\Label{2.17}
\frac{\left|\partial_{z_2}^\nu f(z_1,\delta-\sigma+i\delta)\right|^{\frac1\nu}}{\nu!}\leq e^{\frac\alpha2}\delta^{-1}\sigma^{-\kappa|y_1|}.
\end{equation}
Define $\epsilon=\epsilon_\delta$ by
$$
\epsilon_\delta:=\frac{\T{log}(1-\alpha)+\frac\alpha2}{\kappa\T{log}\,\sigma}.
$$
We note that under the condition $y_1<\epsilon_\delta$ we have $e^{\frac\alpha2}\sigma^{-\kappa y_1}<(1-\alpha)^{-1}$ and hence the term in the right of \eqref{2.17} is $\leq \delta^{-1}(1-\alpha)^{-1}$ that is the radius of convergence of \eqref{2.13} is $\geq \delta':=(1-\alpha)\delta$. 
 The convergence is normal for $z_2$ satisfying $|z_2-((\delta-\sigma)+i\delta)|<\delta'$, uniform in $z_1$ satisfying $x_1\in I,\,\,\eta<y_1<\epsilon$. 
 In fact, by repeating the argument of Theorem~\ref{hartogslemma} in the half-disc $y_1<0$, it is also uniform for $-\epsilon<y_1<-\eta$.
 Hence, by letting $\eta\to 0$,  it remains well defined a holomorphic function on 
$$
(I+i\dot I_\epsilon)\times \Delta_{\delta'}(\delta-\sigma+i\delta)
$$
 that we will denote by $F$. We can of course ``move backwords" and center the Taylor series at $-\delta+\sigma+i\delta$ instead of $\delta-\sigma+i\delta$.
 
 \noindent
 {\bf (g)} In the next step we center the Taylor expansion of $F$ at $\delta+\delta'-2\sigma+i\delta$; by the same argument as in (f) - in which we apply Theorem~\ref{hartogslemma} to a covering of $U^+_\epsilon$ by rescaled discs $x_1+\frac\epsilon\delta \Delta^+_\delta$ - we get extension of $F$ to the domain
 $$
 (I+i\dot I_{\epsilon^2})\times \Delta_{\delta'}(\delta+\delta'-2\sigma+i\delta).
 $$
After finitely many, say N, steps we get $N\delta'-N\sigma>1$; we can also ``move backwords" and get $-N\delta'+N\sigma<-1$. We have thus obtained a holomorphic extension $F$ to a domain of the type
$$
(I+i\dot I_{\epsilon^N})\times \left(\cup_{j=-N}^N\Delta_{\delta'}(j(\delta'-\sigma)+i\delta)\right).
$$
Note that the union of discs in the second term above contains an open neighborhood $V_\delta$ of $I+i\delta$. 
 We move now $\delta$ to $0$: we then get a holomorphic function $F$ defined on the domain
$$
\cup_\delta(I+i\dot I_{\epsilon_\delta^N})\times V_\delta.
$$
We show now that $F$ is also defined for $y_1=0$ that it extends to a wedge $W=I^2+i\Gamma$ where $\Gamma$ is a conic  neighborhood of the axis $y_2>0$ and that it is uniformly continuous at $y=0$ with limit $f$.
We start by remarking that $F$ is defined from the beginning over $U_\delta\times \Delta^+_\delta$ where 
it is uniformly continuous up to $I\times I_\delta$ with limit $f$:
$$
F(z_1,z_2)|_{y_1=0,\,y_2=0^+}=f(x_1,x_2)\quad\T{for  $(x_1,x_2)\in I\times I_\delta$}.
$$
In particular, $F(x_1,z_2)|_{y_2=0^+}=f(x_1,x_2)\quad\T{for $x_2\in I_\delta$}$.
But it then follows from the identity principle for holomorphic functions
\begin{equation*}
F(x_1,z_2)=f(x_1,z_2)\quad\T{for $z_2\in V_\delta\cup \Delta^+_\delta$}.
\end{equation*}
Thus $F$ is defined also for $y_1=0$ and is, wherever defined, an extension of $f(x_1,z_2)$ for values of $z_1$ such that $\pm y_1> 0$. In particular, it is an extension of  $f$ to $\cup_\delta(I+iI_{\epsilon_\delta^N})\times V_\delta$. 
 Next, by a theorem by Kashiwara (cf. \cite{K72}) $F$ extends to a wedge $W$ as before described.
At this point, $F$ has at its own right a  generalized boundary value.
On the other hand, by what we remarked before (d) of the proof of Theorem \ref{2.1}, $F$ is uniformly continuous up to $y=0$ over $I\times D$ where $D$ is open dense in $I$. In particular, for any $x_2^o\in I$ there are neigborhoods 
$B\supset\supset A\ni x^o_2$ 
 such that $ B\setminus A\subset D$. Thus  $F$ is indeed uniformly continuous according to Theorem~\ref{t1.3} and it has limit $f$ at $I^2$. The proof is complete.
\newline
{\bf Proof of Corollary~\ref{c2.1}.}
We let $n=2$ and suppose that $f$ extends to $\{x_1\}\times\Delta^+$ and $\Delta_{\epsilon_{z_2}}\times\{z_2\}$ for any $x_1\in I$ and $z_2\in\Delta^+$ respectively. For $y_2\in I^+$, let $I(y_2)=I+iy_2$; the idea of the proof is to apply Theorem~\ref{2.2} with $I^2$ replaced by $I\times I(y_2)$. But we have to overcome the problem that $f|_{I\times I(y_2)}$ is not known to be continuous. To this end we define $J_l(y_2)=\{z_2\in I(y_2): f(\cdot,z_2)\T{ is holomorphic in $U_{\frac1l}$ and $\underset{z_1\in U_{\frac1l}}{\T{\rm sup}}|f(z_1,z_2)|\leq l$}\}$. We can then see that
$$
f|_{I\times J_l(y_2)}\T{ is continuous}.
$$
If we manage to prove this, then $f$ will be in fact continuous on $U_\delta\times J_l(y_2)$ for some $\delta<\frac1l$ by Theorem~\ref{t1.1} and
this will suffice to carry out our proof. To prove this claim, we write for any pair of points $\tilde x_2+iy_2$, $x_2+iy_2$ in $J_l(y_2)$ and
$\tilde x_1$, $x_1$ in $I$:
\begin{equation*}
\begin{split}
|f(\tilde x_1,\tilde x_2+iy_2)-f(x_1,x_2+iy_2)|&\leq |f(\tilde x_1,\tilde x_2+iy_2)-f(x_1,\tilde x_2+iy_2)|
\\
&+|f(x_1,\tilde x_2+iy_2)-f(x_1,x_2+iy_2)|.
\end{split}
\end{equation*}
Now, the first term is small, for $|\tilde x_1-x_1|$ small and uniformly with respect to $\tilde x_2$, by Cauchy inequalities and by the uniform bound $|f|\leq l$. The second is small for $|\tilde x_2-x_2|$ small  since $f(x_1,\cdot)$ is holomorphic, henceforth continuous.

\end{document}